\documentclass[12pt]{amsart}

\usepackage{amsmath,amsthm,amssymb}
\font\teneufm=eufm10
\font\seveneufm=eufm7
\font\fiveeufm=eufm5
\newfam\frakturfam

\textfont\frakturfam=\teneufm
\scriptfont\frakturfam=\seveneufm
\scriptscriptfont\frakturfam=\fiveeufm

\newtheorem{proposition}{Proposition}

\newtheorem{lemma}{Lemma}
\newtheorem{theorem}{Theorem}
\newtheorem{corollary}{Corollary}
\newtheorem{remark}{Remark}

\def\bee{\begin{eqnarray}}
\def\bes{\begin{eqnarray*}}
\def\eee{\end{eqnarray}}
\def\ees{\end{eqnarray*}}

\def\s{\sigma}
\def\t{\tau}

\def\Proof{{\sl Proof.}\ }

\title{Wild Automorphisms of Free Metabelian Lie and Associative Algebras}

\begin{document}
\date{}

\maketitle

\begin{center}

{\bf Ualbai Umirbaev}\footnote{Department of Mathematics,
 Wayne State University,
Detroit, MI 48202, USA; Center for Theoretical and Applied Algebra, 
Al-Farabi Kazakh National University, Almaty, 050040, Kazakhstan,
e-mail: {\em umirbaev@wayne.edu}}

\end{center}

\begin{abstract} We show that the Chein automorphism $\tau=(x_1+[[x_2,x_3],x_1],x_2,\ldots,x_n)$ of the free metabelian Lie algebra of rank $n\geq 3$ over an arbitrary field $K$ is wild. We also prove that the Chein automorphisms $\s=(x_1+[x_2,x_3]x_1,x_2,\ldots,x_n)$ and $\tau$ of the free metabelian associative algebra of rank $n\geq 3$ over $K$ are wild. 
\end{abstract}



\tableofcontents

\section{Introduction}

\hspace*{\parindent}

Let $L$ be the free metabelian Lie algebra of rank $n$ over a field $K$, with free generators $x_1,x_2,\ldots,x_n$. 
As usual, we write $\phi=(f_1,\ldots,f_n)$ for the endomorphism of $L$ defined by $\phi(x_i)=f_i$ for all $1\leq i\leq n$. 
 In 2015, Nauryzbaev \cite{Nauryzbaev15} proved that if $n\geq 4$ and $K$ has characteristic different from $3$, then the tame automorphism group $\mathrm{TAut}(L)$ of $L$ is generated by all linear automorphisms together with the quadratic automorphism
\bes
\epsilon=(x_1+[x_2,x_3],x_2,\ldots,x_n). 
\ees
If $n\geq 5$, the same result holds over fields of arbitrary characteristic.

It is shown in \cite{AutoI,AutoII} that if $n\geq 4$ and $K$ has characteristic different from $3$, then the full group of automorphisms $\mathrm{Aut}(L)$ of $L$ is generated by all linear automorphisms together with the quadratic automorphism $\epsilon$ 
and the cubic Chein automorphism
\bes
\tau=(x_1+[[x_2,x_3],x_1],x_2,\ldots,x_n). 
\ees
If $n\geq 5$, the same result holds over fields of arbitrary characteristic.  

It is well known that $\tau$ is not tame when $n=3$ \cite{Nauryzbaev09}. It can be shown that $\tau$ is {\em absolutely} wild when $n=3$ and $K$ is a field of characteristic zero \cite{Drensky92,KP16,SU25}.

For $n\geq 4$, if $K$ is a field of characteristic zero, then $\tau$ is {\em approximately tame}, since every automorphism of $M_n$ is approximately tame in this case \cite{BD93-2,KP16,SU25}.

The main result of this paper establishes that the Chein automorphism $\tau$ is wild over any field $K$. Consequently, $\epsilon$ and $\tau$ form a minimal generating set for $\mathrm{Aut}(M_n)$ modulo linear transformations. 

We first establish the wildness of certain automorphisms of free metabelian associative algebras. Recall that the variety of associative algebras defined by the identity
\bes
[x,y][z,t]=0
\ees
is called {\em metabelian} \cite{95SMJ}. 

Let $M$ be the free metabelian associative algebra of rank $n$ over $K$ with free generators $x_1,\ldots,x_n$.  It is not difficult to show that, when $n=2$, every nontrivial $\mathrm{IA}$-automorphism of $M$ (i.e., an automorphism inducing the identity modulo the commutator ideal) is wild. The well-known Anick automorphism of the free associative algebra (see \cite{UU07JR}) induces a wild automorphism of $M$ when $n=3$ and $K$ is a field of characteristic zero. Further examples of wild automorphisms in rank three are given in \cite{DY,Liu}. To the best of our knowledge, the case $n\geq 4$ has not previously been studied.

The free metabelian Lie algebra $L$ can be identified with the Lie subalgebra of the commutator algebra $M^{(-)}$ generated by  
$x_1,\ldots,x_n$. Thus $\tau$ is also an automorphism of $M$. We prove that the Chein automorphisms 
\bes
\s=(x_1+[x_2,x_3]x_1,x_2,\ldots,x_n)
\ees
 and $\tau$ of the free metabelian associative algebra $M$ of rank $n\geq 3$ over $K$ are both wild.

The paper is organized as follows. In Section 2, we establish a criterion for the invertibility of matrices of a certain type over free associative algebras. In Section 3, we show that certain endomorphisms of the free associative algebra cannot be automorphisms by proving that their components do not generate the augmentation ideal as a left module. In Section 4 we prove that the automorphisms $\sigma$ and $\tau$ of the free metabelian associative algebra $M$ are wild. In Section 5, we state the wildness of the automorphism $\tau$ of $L$ as a separate theorem, since this result completes a series of investigations by the author into automorphisms of free metabelian Lie algebras (see \cite{AutoI,AutoII}).

\section{Invertibility of Matrices}

\hspace*{\parindent} 

Let $A=K\langle y_1,\ldots,y_n\rangle$ be the free associative algebra over a field $K$ in the variables $y_1,\ldots,y_n$. Let 
\bes
T=I_m+T_p+T_{p+1}+\cdots+T_q, 1\leq p\leq q, 
\ees
be an $m\times m$ matrix over $A$, where $I_m$ is the identity matrix, and $T_i$ is a matrix of homogeneous elements of degree $i$ for all $p\leq i\leq q$. Bryant and Drensky proved that if $T$ is invertible, then $\mathrm{Tr}(T_p)\in [A,A]$, where $\mathrm{Tr}$ denotes the trace \cite[Proposition 3.7]{BD93-1} (see also \cite[Lemma 3.4]{Liu} or \cite[Lemma 7]{SU25}). 

We are going to prove one variation of this result. 
 Let $C$ be the commutator ideal of $A$. 

\begin{proposition}\label{p1}
	Let 
	\bes
	S=I_m+T, 
	\ees
	be an $m\times m$ matrix, where all entries of $T$ lie in $C$.  If $S$ is invertible, then 
\bes
\mathrm{Tr}(T)\in [A,C]+C^2. 
\ees
\end{proposition}

Let $R$ be an arbitrary unital augmented associative $K$-algebra and let $\eta : R\to K$ be its augmentation map. Denote by $M_m(R)$ the set of all $m\times m$ matrices over $R$. Let $GL_m(R)$ be the group of all invertible matrices from $M_m(R)$. Then $SL_m(K)$ is the kernel of the homomorphism $GL_m(R)\to GL_m(K)$ induced by $\eta$. If $J$ is an ideal of $R$ then relative general linear group $GL_m(R,J)$ and relative special linear group $SL_m(R,J)$ are defined as kernels of the homomorphisms 
\bes
GL_m(R)\to GL_m(R/J)
\ees
 and 
 \bes
 SL_m(R)\to SL_m(R/J),
 \ees
 respectively, induced by the canonical projection $R\to R/J$.

Let $E_{ij}$ denote the $m\times m$ matrix with $1$ in the $(i,j)$-entry and $0$ elsewhere. 
A matrix of the form 
\bes
e_{ij}(a)=I+aE_{ij}, \ i\neq j, \ a\in R,
\ees
is called {\em elementary}. Denote by $E_m(R)$ the subgroup of $SL_m(R)$ generated by all elementary matrices. 
If $J$ is an ideal of $R$ then denote by $E_m(R,J)$ the normal closure in $E_m(R)$ of the set of all elementary matrices of the form $e_{ij}(a)$, where $a\in J$. 

\begin{lemma}\label{l1}  $SL_m(A)=E_m(A)$. 
\end{lemma}
\Proof Every invertible matrix over a free associative algebra can be expressed as a product of elementary matrices and a diagonal matrix of the form $\mathrm{diag}=(d,1,\ldots,1)$, where $d$ is the determinant of the constant part of the original matrix \cite[Theorem 2.4.4]{Cohn}. If $M\in SL_m(A)$, then the determinat of its constant part $M(0)$ is $1$, i.e., $d=1$ in this case. $\Box$

The factor algebra $B=A/C=K[z_1,\ldots,z_n]$ is the polynomial algebra in the variables $z_1=y_1+C,\ldots,z_n=y_n+C$. Let 
\bes
\pi: A\to B
\ees
be the natural projection. This homomorphism induces the homomorphism 
\bes
\pi^*: SL_m(A)\to SL_m(B). 
\ees

In 1977, Suslin \cite{Suslin} proved  that $SL_m(B)=E_m(B)$ for all $m\geq 3$. 

Recall that the elementary matrices over any associative unital algebra $R$ satisfy  the well known \emph{Steinberg relations} (see \cite{Milnor}): 

$(S1)$ $e_{ij}(a)e_{ij}(b)=e_{ij}(a+b) \quad (i\ne j)$;

$(S2)$ $[e_{ij}(a),e_{\ell r}(b)]=1 \quad 
	(j\ne\ell,\ i\ne r)$,\\
	where $[g,h]=ghg^{-1}h^{-1}$; 
	
$(S3)$ $[e_{ij}(a),e_{j\ell}(b)]=e_{i\ell}(ab) \quad 
	(i,j,\ell\ \text{pairwise distinct})$.

Elementary matrices over fields also satisfy the following relation (see \cite{Milnor}): 

$(D)$ $\{u,v\}=I_m \quad 0\neq u,v\in K$, \\
where 
\bes
\{u,v\}=h_{ij}(uv) h_{ij}(u)^{-1} h_{ij}(v)^{-1}, \\
h_{ij}(u)=w_{ij}(u) w_{ij}(-1), \\
w_{ij}(u)=X_{ij}(u) X_{ji}(-u^{-1}) X_{ij}(u). 
\ees

In 1982, Tulenbaev \cite{Tulenbaev} proved that, for $m\geq 5$, the relations $(S1)$--$(S3)$ and $(D)$ form a complete set of defining relations for the elementary generators of $SL_m(B)=E_m(B)$.

\begin{lemma}\label{l2} 
	$\mathrm{Ker}(\pi^*)=E_m(A,C)$ for all $m\geq 5$. 
\end{lemma}
\Proof We have $e_{ij}(c)\in \mathrm{Ker}(\pi^*)$ if $c\in C$. This implies that $E_m(A,C)\subseteq \mathrm{Ker}(\pi^*)$. Then $\pi^*$ induces the homomorphism 
\bes
\lambda : SL_m(A)/E_m(A,C)=E_m(A)/E_m(A,C)\to SL_m(B)=E_m(B). 
\ees
Let $s:B\to A$ be a section of $\pi$, that is, $\pi\circ s=\mathrm{Id}_B$. Every $a\in A$ can be uniquely written in the form 
\bes
a= (s\circ\pi)(a) +c, c\in C. 
\ees
Define a homomorphism 
\bes
\rho : SL_m(B)=E_m(B)\to E_m(A)/E_m(A,C)=SL_m(A)/E_m(A,C), 
\ees
by 
\bes
\rho(e_{ij}(\pi(a)))= e_{ij}((s\circ\pi)(a))E_m(A,C). 
\ees
It is possible since the elementary matrices over $A$ satisfy the relations $(S1)$--$(S3)$ and $(D)$ and the same relations are defining relations for $SL_m(B)=E_m(B)$. 

It is clear that $\lambda$ and $\mu$ are mutually inverse. Consequently, $\lambda$ is an isomorphism and $\mathrm{Ker}(\pi^*)=E_m(A,C)$. $\Box$

\begin{remark} The proof of Lemma \ref{l2} follows that of Theorem 3.1 in \cite{UU07JR}. Moreover, the statement remains valid for any unital algebra $R$ equipped with a surjective homomorphism $\pi :R\to B$ inducing the homomorphism $\pi^*: E_m(R)\to SL_m(B)=E_m(B)$. 
\end{remark}

\begin{lemma}\label{l3} Let $M\in M_m(C)$ be a matrix  such that $\mathrm{Tr}(M)\in [A,C]$. 
Then, for every elementary matrix $e_{ij}(a)\in E_m(A)$, we have
\bes 
\mathrm{Tr}(e_{ij}(a)Me_{ij}^{-1}(a))\in [A,C]. 
\ees
\end{lemma}
\Proof 
Let $M=[c_{kl}]$. Observe that $e_{ij}^{-1}(a)=I_m-aE_{ij}$. Then 
\bes
e_{ij}(a)Me_{ij}^{-1}(a)=M+ aE_{ij}M-MaE_{ij}-aE_{ij}MaE_{ij}.  
\ees
We have 
\bes
aE_{ij}M=\sum_{r}ac_{jr}E_{ir}, \quad MaE_{ij}=\sum_{r}c_{ri}aE_{rj}, \quad aE_{ij}MaE_{ij}=ac_{ji}aE_{ij}. 
\ees

Therefore 
\bes
\mathrm{Tr}(aE_{ij}M)=ac_{ji},\quad 
\mathrm{Tr}(MaE_{ij})=c_{ji}a,\quad 
\mathrm{Tr}(aE_{ij}MaE_{ij})=0. 
\ees

Consequently, 
\bes
\mathrm{Tr}(e_{ij}(a)Me_{ij}^{-1}(a))=\mathrm{Tr}(M)+ac_{ji}-c_{ji}a=\mathrm{Tr}(M)+[a,c_{ji}]\in [A,C]. \quad \Box
\ees

{\bf Proof of Proposition \ref{p1}}. Let $S=I_m+T$ be a matrix satisfying the hypothesis of Proposition \ref{p1}, that is, $T\in M_m(C)$. We may assume that $m\geq 5$, since enlarging the matrix $T$ by adjoining zero rows and columns does not change its trace. If $S$ is invertible, then $S$ belongs to the kernel of $\pi^*$. By Lemma \ref{l2}, 
\bes
S=S_1\ldots S_k, 
\ees 
where each $S_i$ is a matrix of the form 
\bes
S_i= Q_i e_{p_iq_i}(c_i) Q_i^{-1}=I_m+Q_ic_iE_{p_iq_i} Q_i^{-1}
\ees
for some $1\leq p_i\neq q_i\leq m$, $c_i\in C$, and $Q_i\in E_m(A)$. 

We have $\mathrm{Tr}(c_iE_{p_iq_i})=0$. Since $Q_i\in E_m(A)$ is a product of elementary matrices, repeated application of Lemma \ref{l3} shows that the trace of the matrix 
\bes
T_i=Q_ic_iE_{p_iq_i} Q_i^{-1}
\ees
belongs to $[A,C]$. 

Further, 
\bes
I_m+T=(I_m+T_1)\ldots (I_m+T_k)
\ees
implies that 
\bes
T=T_1+\ldots +T_k +M, 
\ees
where
\bes
M=\sum_{i_1<\ldots<i_s,s\geq 2}T_{i_1}\ldots T_{i_s}
\ees
Then $\mathrm{Tr}(M)\in C^2$ since $M\in M_m(C^2)$. Consequently, 
\bes
\mathrm{Tr}(T)=\mathrm{Tr}(T_1)+\ldots +\mathrm{Tr}(T_k) +\mathrm{Tr}(M)\in [A,C]+C^2. \quad \Box
\ees

\section{Nongenerators of the Augmentation Ideal}

\hspace*{\parindent} 

Denote by $I_A$ the augmentation ideal of $A$. Then 
\bes
I_A=Ay_1\oplus\ldots\oplus Ay_n
\ees
is a free left $A$-module. Set $Y=(y_1,\ldots,y_n)$. For any $a\in I_A$ there exists a unique vector $v\in A^n$ such that $a=vY^t$, where $t$-means the transposition. 

\begin{corollary}\label{c1} Let 
\bee\label{f1}
S=I_n + D+T, 
\eee	
where 
\bes
D=\mathrm{diag}([y_2,y_3],0,\ldots,0),  
\ees
and $T\in M_n(J)$ with $J=I_AC+CI_A$. Then $S$ is not invertible.  
\end{corollary}
\Proof We have  $D+T\in M_n(C)$ and 
\bes
\mathrm{Tr}(D+T)=\mathrm{Tr}(D)+\mathrm{Tr}(T)\in [y_2,y_3]+J. 
\ees
Therefore, $\mathrm{Tr}(D+T)\notin [A,C]$. By Proposition \ref{p1}, $S$ is not invertible. $\Box$

\begin{corollary}\label{c2}
$(a)$ The elements
\bes
y_1+[y_2,y_3]y_1 +t_1, y_2+t_2,\ldots,y_n+t_n, 
\ees
of $A$, where $t_i\in C^2$ for all $i$, do not generate $I_A$ as a left $A$-module. 

$(b)$ The  elements 
\bes
y_1+[[y_2,y_3],y_1] +t_1, y_2+t_2,\ldots,y_n+t_n,  
\ees
of $A$, where $t_i\in C^2$ for all $i$, do not generate $I_A$ as a left $A$-module. 
\end{corollary}
\Proof Let 
\bee\label{f2}
\phi=(y_1+[y_2,y_3]y_1 +t_1, y_2+t_2,\ldots,y_n+t_n)
\eee
and 
\bee\label{f3}
\psi=(y_1+[[y_2,y_3],y_1] +t_1, y_2+t_2,\ldots,y_n+t_n). 
\eee

Every element $t_i$ can be uniquely written in the form 
\bes
t_i=t_{i1}x_1+\ldots+t_{in}x_n. 
\ees
Since $t_i\in C^2$ it follows that $t_{ij}\in J$. Set $T=[t_{ij}]$. Then 
\bes
\phi^t=SY^t, 
\ees
where $S$ is the matrix (\ref{f1}). 

If the statement $(a)$ is false, there exists a matrix $P$ such that  
\bes
Y^t=P\phi^t=PSY^t. 
\ees
Since $I_A$ is a free module with basis $y_1,\ldots,y_n$ it follows that $PS=I_n$, i.e., $S$ is invertible. This contradicts Proposition \ref{p1}. 

We have 
\bes
[[y_2,y_3],y_1]=[y_2,y_3]y_1-y_1y_2y_3+y_1y_3y_2. 
\ees
Set 
\bes
a=y_1+[[y_2,y_3],y_1]+t_1+y_1y_2(y_3+t_3)-y_1y_3(y_2+t_2)\\
=y_1+[y_2,y_3]y_1 + (t_1+y_1y_2t_3-y_1y_3t_2)=y_1+[y_2,y_3]y_1+t_1',
\ees
where $t_1'=t_1+y_1y_2t_3-y_1y_3t_2\in C^2$. Observe that the elements given in $(b)$ generate  $I_A$ as a left $A$-module if and only if the elements 
\bes
x_1+[x_2,x_3]x_1+u_1',x_2+u_2,\ldots,x_n+u_n  
\ees
do so. However, this is impossible by part $(a)$. $\Box$

\begin{corollary}\label{c3}
Neither of the endomorphisms $\phi$ and $\psi$, defined by (\ref{f2}) and (\ref{f3}), respectively, is not an automorphism  of $A$. 
\end{corollary}
\Proof If $\phi$ (or $\psi$) is an automorphism, then its components generate $I_A$ as a left $A$-module. This is impossible by Corollary \ref{c2}. $\Box$

\section{Wild Automorphisms of Free Metabelian Associative Algebras}

\hspace*{\parindent} 

The algebra $M=A/C^2$ is the free metabelian associative algebra in the varibles $x_i=y_i+C^2$, where $1\leq i\leq n$ \cite{95SMJ}. 

\begin{theorem}\label{t1} The Chein automorphisms 
\bes
\s=(x_1+[x_2,x_3]x_1,x_2,\ldots,x_n)
\ees
and 
\bes
\t=(x_1+[[x_2,x_3],x_1],x_2,\ldots,x_n)
\ees
of the free metabelian associative algebra $M$ of rank $n\geq 3$ are both wild. 
\end{theorem}
\Proof Suppose that $\s$ (or $\tau$) is tame. Then it can be expressed as a product of elementary automorphisms.  By lifting each elementary factor to an elementary automorphism of $A$, we obtain an automorphism $\phi$ (or $\psi$) of $A$ that induces $\s$ (or $\tau$). Clearly, $\phi$ (or $\psi$) must be of the form (\ref{f2}) (or (\ref{f3}), respectively). However, by Corollary \ref{c3}, neither $\phi$ nor $\psi$ can be automorphisms of $A$, a contradiction. $\Box$

\section{A Wild Automorphism of Free Metabelian Lie Algebras}

\hspace*{\parindent} 

The Lie subalgebra $L$ of the commutator algebra $M^{(-)}$ generated by  $x_1,\ldots,x_n$ is a free metabelian Lie algebra on these generators (see \cite{95SMJ}).

As mentioned in the introduction, if $n\geq 4$ and $K$ has characteristic different from $3$, then the full group of automorphisms $\mathrm{Aut}(L)$ is generated by all linear automorphisms together with the quadratic automorphism $\epsilon$ 
and the cubic Chein automorphism $\tau$ \cite{AutoI,AutoII}. 
If $n\geq 5$, the same result holds over fields of arbitrary characteristic.

The following corollary of Theorem \ref{t1} completes the author's long-standing study of automorphisms of free metabelian Lie algebras (see \cite{AutoI,AutoII}). In view of its significance, we state it as a separate theorem.

\begin{theorem}\label{t2} The Chein automorphism 
	\bes
	\t=(x_1+[[x_2,x_3],x_1],x_2,\ldots,x_n)
	\ees
	of the free metabelian Lie algebra $L$ of rank $n\geq 3$ is wild. 
\end{theorem}
\Proof Since every elementary automorphism of $L$ is also an elementary automorphism of $M$, Theorem \ref{t1} implies that $\t$ cannot be expressed as a product of elementary automorphisms of $L$. $\Box$

\begin{theorem}\label{t3} The automorphisms $\epsilon$ and $\tau$ form a minimal generating set for $\mathrm{Aut}(L)$ modulo linear automorphisms  whenever either $n\geq 4$ and $\mathrm{char}(K)\neq 3$, or $n\geq 5$. 
\end{theorem}
\Proof It is shown in \cite{AutoI,AutoII} that $\epsilon$ and $\tau$ generate $\mathrm{Aut}(L)$ modulo linear automorphisms. Theorem \ref{t2} establishes that $\tau$ does not belong to the subgroup generated by $\epsilon$ and all linear automorphisms. 

Observe that $\epsilon$ is quadratic, whereas $\tau$ is cubic (see definitions in \cite{AutoI,AutoII}). Every $\mathrm{IA}$ -automorphism belonging to the subgroup generated by $\tau$ and all linear automorphisms can be expressed as a product of linear conjugates of $\tau$ and its inverse 
\bes
\t^{-1}=(x_1-[[x_2,x_3],x_1],x_2,\ldots,x_n)
\ees
 (see the proof of Lemma 1 in \cite{AutoII}). Consequently, this subgroup contains no quadratic $\mathrm{IA}$-automorphisms and, in particular, does not contain $\epsilon$. $\Box$

\section*{Acknowledgments}

AI tools were used to improve the language of the manuscript.

\end{document}